\documentclass[12pt]{article}

\usepackage{acronym}
\usepackage{amsfonts}
\usepackage{amsmath,amssymb,bm}
\usepackage{graphicx}
\usepackage{subfigure}
\usepackage{mdframed}
\usepackage{url}
\usepackage{multirow}

\newtheorem{theorem}{Theorem} 

\newtheorem{proposition}{Proposition}[section]
\newtheorem{lemma}{Lemma}[section]
\newtheorem{definition}{Definition}

\def\be{\begin{Example}}
\def\ee{\end{Example}}
\def\bt{\begin{theorem}}
\def\et{\end{theorem}\bigskip}
\def\bl{\begin{Lemma}}
\def\el{\end{Lemma}\bigskip}
\def\ep{\end{Proposition}\bigskip}
\def\bp{\begin{Proposition}}
\def\bd{\begin{definition}}
\def\ed{\end{definition}}

\begin{document}
\title{\bf Spectral projected gradient method for generalized tensor eigenvalue complementarity problem
 }
 \author{
 Gaohang Yu \thanks{\small
  School of Mathematics and Computer Sciences,
Gannan Normal University, Ganzhou, 341000, China. E-mail:
maghyu@163.com}
\quad Yisheng Song \thanks{\small
 School of Mathematics and Information Science and Henan Engineering Laboratory for Big Data Statistical Analysis and Optimal Control, Henan Normal University, China. E-mail: songyisheng1@gmail.com}
\quad Yi Xu \thanks{\small
 Department of Mathematics, Southeast University, China. E-mail: yi.xu1983@gmail.com}
 \quad Zefeng Yu \footnotemark[1]
  }

\date{\today}
\maketitle


\begin{abstract}
  This paper looks at the tensor eigenvalue complementarity problem (TEiCP) which arises from the stability analysis of finite dimensional mechanical systems and is closely related to the optimality conditions for polynomial optimization. We investigate two monotone ascent spectral projected gradient (SPG) methods for TEiCP. We also present a shifted scaling-and-projection algorithm (SPA), which is a great improvement of the original SPA method proposed by Ling, He and Qi [Comput. Optim. Appl., DOI 10.1007/s10589-015-9767-z]. Numerical comparisons with some existed gradient methods in the literature are reported to illustrate the efficiency of the proposed methods.


\end{abstract}

\section{Introduction}

A $m$th-order $n$-dimensional real tensor $\mathcal{A}$ consists of
$n^m$ entries in real numbers:
\[
\mathcal{A}=(a_{i_1i_2\cdots
i_m}),\,\, a_{i_1i_2\cdots i_m}\in \mathbb{R},\,\, \mbox{for any}  \ \ i_1,i_2,\ldots,i_m \in [n],
\]
where $[n]=\{1,2,\ldots,n\}.$ Denote the set of all real $m$th-order $n$-dimensional tensors by $\mathbb{T}^{[m,n]}$. $\mathcal{A}$  is called \textit{symmetric} if the value of
$a_{i_1i_2\cdots i_m}$ is invariant under any
permutation of its indices $i_1,i_2,\ldots,i_m$. Denote the set of all real symmetric $m$th-order $n$-dimensional tensors by $\mathbb{S}^{[m,n]}$.
For any vector $x\in \mathbb{R}^n$, $\mathcal{A}x^{m-1}$ is a vector in $\mathbb{R}^n$ with its $i$th
component as
$$
(\mathcal{A}x^{m-1})_i=\sum_{i_2,\ldots,i_m=1}^n a_{i i_2 \cdots i_m}x_{i_2}\cdots
x_{i_m}.
$$
 A real symmetric tensor $\mathcal{A}$ of order $m$ dimension $n$ uniquely defines a
$m$th degree homogeneous polynomial function $h$ with  real
coefficient by
$$
h(x):=\mathcal{A}x^m=x^T(\mathcal{A}x^{m-1})=\sum_{i_1,\ldots, i_m=1}^n a_{i_1 \cdots i_m}x_{i_1}\cdots
x_{i_m}.
$$

We call that the tensor $\mathcal{A}$ is positive definite if $\mathcal{A}x^m>0$ for all $x\neq 0$.

In 2005, Qi \cite{Qi05} and Lim \cite{Lim05} proposed the definition of eigenvalues and eigenvectors for higher order tenors, independently.
Furthermore, in \cite{CPZ09}, these definitions were unified by Chang, Person and Zhang. Let $\mathcal{A}$ and $\mathcal{B}$ be real-valued, $m$th-order $n$-dimensional symmetric tensors. Assume further that $m$ is even and $\mathcal{B}$ is positive definite. we call $(\lambda, x) \in \mathbb{R} \times \mathbb{R}^n \backslash \{0\}$ is a {\bf generalized eigenpair} of $(\mathcal{A},\mathcal{B})$ if
\begin{equation}\label{generalized eigenpair}
\mathcal{A}x^{m-1}=\lambda\mathcal{B}x^{m-1}.
\end{equation}

    When the tensor $\mathcal{B}$ is an identity tensor $\mathcal{\varepsilon}$ such that $\mathcal{\varepsilon}x^{m-1}=\|x\|^{m-2}x$ for all $x \in \mathbb{R}^n$ \cite{CPZ09}, the eigenpair reduces to $Z$-eigenpair \cite{Qi05,Lim05}. 
 Another special case is that when $\mathcal{B}=\mathcal{I}$ with
 $$
 (\mathcal{I})_{i_1 i_2\cdots i_m}=\delta_{i_1 i_2\cdots i_m}=\left\{\begin{array}{ll}
1,& \ \mbox{if} \ i_1=i_2=\ldots=i_m,\\
0,& \mbox{otherwise},\\
\end{array}\right.
 $$
 the real scalar $\lambda$ is called an $H$-eigenvalue and the real vector $x$ is the associated $H$-eigenvector of the tensor $\mathcal{A}$ \cite{Qi05}. In the last decade, tensor eigenproblem has received much attention in the literature \cite{CDN14,DW15,HCD15,HCD15-PJO,LQY13,LQY13-LAA,SQ15,SY15,YYZ14,YY10}, which has numerous applications \cite{CDHS13,QWW08,QYW,QYX}.

In this paper, we consider the tensor eigenvalue complementarity problem (TEiCP):

\begin{mdframed}
finding a scalar $\lambda\in \mathbb{R},$ and $x\in \mathbb{R}^n \backslash \{0\}$ such that\\
\begin{equation}\label{TEiCP.problem}
\begin{split}
x\ge0,\\
(\lambda\mathcal{B}-\mathcal{A})x^{m-1}\ge0, \\
\langle x,(\lambda\mathcal{B}-\mathcal{A})x^{m-1}\rangle =0, \\
\end{split}
\end{equation}
where $\mathcal{A}\in \mathbb{T}^{[m,n]}$, and $\mathcal{B}\in \mathbb{S}^{[m,n]}$ is positive definite.
\end{mdframed}
 The solution of TEiCP $(\lambda,x)$ is called Pareto eigenpair of $(\mathcal{A},\mathcal{B})$. In some special case, we can call it Pareto H-eigenpair or Pareto Z-eigenpair  \cite{SQ15-JGO} if the tensor $\mathcal{B}$ has special form as shown above in the generalized eigenpairs (\ref{generalized eigenpair}). Replacing the nonnegative cones in (\ref{TEiCP.problem}) by a closed convex cone and its dual cone,  Ling, He and Qi investigated the cone eigenvalue complementarity problem for higher-order tensor in \cite{LHQ2015-COAP1}. Moreover, in \cite{LHQ2015-COAP2}, they studied the high-degree eigenvalue complementarity problem for tensors as a natural extension of quadratic eigenvalue complementarity problem for matrices. TEiCP is also closely related to the optimality conditions for polynomial optimization \cite{SQ15-JGO}, a class of differential inclusions with noncovex processes \cite{LHQ2015-COAP1}, and a kind of nonlinear differential dynamical system \cite{ChenYangYe15}.
The properties of Pareto eigenvalues and their connection to polynomial optimization are studied in \cite{SQ15-JGO}. Recently, as a special type of nonlinear complementarity problems, the tensor complementarity problem is inspiring more and more research in the literature \cite{BHW15,CQW15,Z.ChenQi15,ChenYangYe15,DLQ15,GLQX15,LQX15,SQ15,SY15,WHB15}. A shifted projected power method for TEiCP was proposed in \cite{ChenYangYe15}, in which they need an adaptive shift to force the objective to be (locally) convex to guarantee the convergence of power method. In \cite{LHQ2015-COAP1}, Ling, He and Qi presented a scaling-and-projection algorithm (SPA) for TEiCP. One main shortcoming of SPA is the stepsize will approach to zero as the sequence gets close to a solution of TEiCP \cite{LHQ2015-COAP1}. Recently, by introducing an NCP-function, Chen and Qi \cite{Z.ChenQi15} reformulated the TEiCP as a system of nonlinear equations. And then, they proposed a semismooth Newton method for solving the system of nonlinear equations \cite{Z.ChenQi15}.

 In this paper, we will investigate two spectral projected gradient algorithms for TEiCP.  The rest of this paper is organized as follows. In Section 2, some properties of the solutions of TEiCP and two optimization reformulations of TEiCP are presented. In Section 3, two spectral projected gradient algorithms are proposed. Global convergence results could be established under some suitable assumptions. We also present a shifted scaling-and-projection algorithm (SSPA) in Section 4, which is a great improvement of the original SPA method \cite{LHQ2015-COAP1}. Numerical experiments are reported in Section 4 to show the efficiency of the proposed methods. Finally, we have a conclusion section.

 Throughout this paper, let $\mathbb{R}^n_+=\{x\in \mathbb{R}^n:x\ge 0\}$, and $\mathbb{R}^n_{++}=\{x\in \mathbb{R}^n:x> 0\}$. Given a set $J\subseteq [n]$, the principal sub-tensor of a tensor $\mathcal{A}\in \mathbb{T}^{[m,n]}$, denoted by $\mathcal{A}_J$, is tensor in $\mathbb{T}^{[m,|J|]}$, such that $\mathcal{A}_J=(a_{i_1\ldots i_m})$ for all $i_1,\ldots, i_m \in J$. Here, the symbol $|J|$ denotes the cardinality of $J$.

\section{Some properties and reformulations of TEiCP}
The following proposition shows the relationship between the solution of TEiCP (\ref{TEiCP.problem}) and the generalized eigenvalue problem (\ref{generalized eigenpair}).
\begin{proposition}
$(\lambda,x)$ is a solution of TEiCP (\ref{TEiCP.problem}) if and only if there exists a subset $I\subseteq[n]$, such that $\lambda$ is a generalized eigenvalue of $(\mathcal{A}_I, \mathcal{B}_I)$ and $x_I\in \mathbb{R}^{|I|}_{++}$ is a corresponding eigenvector, and
$$\sum_{i_2,\ldots,i_m\in I}(\lambda b_{ii_2\cdots i_m}-a_{ii_2\cdots i_m})x_{i_2}\cdots x_{i_m}\ge 0, \; \forall \; i \in \bar{I}:=[n]\backslash I.$$
In such a case, the Pareto eigenvector $x$ satisfies $x_{\bar{I}}=0$.
\end{proposition}
This proposition was firstly presented in \cite{SQ15-JGO} for Pareto H-eigenpair and Pareto Z-eigenpair, and then unified by Xu and Ling for TEiCP \cite{XL15-ORC}.

Denote the set of solutions of (\ref{TEiCP.problem}) by $\sigma (\mathcal{A}, \mathcal{B})$ , i.e.,
$$\sigma (\mathcal{A}, \mathcal{B})=\{(\lambda,x)\in \mathbb{R}\times\mathbb{R}^n \backslash \{0\}:0\le x \bot (\lambda\mathcal{B}-\mathcal{A})x^{m-1}\ge0 \}.$$

 If $(\lambda,x)\in \sigma (\mathcal{A}, \mathcal{B})$, then $(\lambda,sx)\in \sigma (\mathcal{A}, \mathcal{B})$ for any $s>0$. On the other hand, given a tensor $\mathcal{A}\in \mathbb{T}^{[m,n]}$,  we know that there exists the unique semi-symmetric tensor \cite{NQ15-GlobalOpt} $\bar{\mathcal{A}}$ such that $\mathcal{A}x^{m-1}=\bar{\mathcal{A}}x^{m-1}$. It is clear that $\sigma (\mathcal{A}, \mathcal{B})=\sigma (\bar{\mathcal{A}}, \mathcal{B}) $. Without loss of generality, we always assume that $\mathcal{A}\in \mathbb{S}^{[m,n]}$ and just consider the solutions on the unit-sphere with $\|x\|_2=1$.

 \begin{proposition}
 The symmetric TEiCP (\ref{TEiCP.problem}) is equivalent to the following optimization problem 
\begin{equation}\label{max-optimization-problem}
\max \lambda(x)=\frac{\mathcal{A}x^m}{\mathcal{B}x^m} \ \ \mbox{subject to} \ \ x \in \mathbb{S}^{n-1}_+:=\{x\in\mathbb{R}^n:x^Tx=1, x\ge0\},
\end{equation}
in the sense that any equilibrium solution $x$ of (\ref{max-optimization-problem}) is a solution of the symmetric TEiCP.
\end{proposition}

By some simple calculations, we can get its gradient and Hessian are as follows
\begin{equation}\label{gradient}
g(x)\equiv \nabla \lambda(x)=\frac{m}{\mathcal{B}x^m}(\mathcal{A}x^{m-1}-\frac{\mathcal{A}x^m}{\mathcal{B}x^m}\mathcal{B}x^{m-1}).
\end{equation}
and its Hessian is
\begin{align*}\label{Hessian}
H(x)\equiv & \nabla^2 \lambda(x)\\
&=\frac{m(m-1)\mathcal{A}x^{m-2}}{\mathcal{B}x^m}-\frac{m(m-1)\mathcal{A}x^m\mathcal{B}x^{m-2}+m^2(\mathcal{A}x^{m-1}\circledcirc\mathcal{B}x^{m-1})}{(\mathcal{B}x^m)^2}\\
&+\frac{m^2\mathcal{A}x^m(\mathcal{B}x^{m-1}\circledcirc\mathcal{B}x^{m-1})}{(\mathcal{B}x^m)^3},
\end{align*}
where $x\circledcirc y= xy^T+yx^T$, and $\mathcal{A}x^{m-2}$ is a matrix with its component as $$(\mathcal{A}x^{m-2})_{i_1i_2}=\sum_{i_3,\ldots,i_m=1}^n a_{i_1 i_2i_3 \cdots i_m}x_{i_3}\cdots
x_{i_m}\ \ \mbox{for all}\ \ i_1,i_2\in [n].$$

According to (\ref{gradient}), we can derive that the gradient $g(x)$ is located in the tangent plane of  $\mathbb{S}^{n-1}$ at $x$, since
 \begin{equation}\label{xtgx}
x^Tg(x)=\frac{m}{\mathcal{B}x^m}(x^T\mathcal{A}x^{m-1}-\frac{\mathcal{A}x^m}{\mathcal{B}x^m}x^T\mathcal{B}x^{m-1})=0.
\end{equation}

The Lagrangian function is $L(x,\mu,v)=\lambda(x)+\mu (x^Tx-1)+v^Tx,$ where $\mu  \in \mathbb{R}$ and $v \in \mathbb{R}^n$ are the Lagrange multipliers. Any equilibrium solution of the nonlinear programming problem (\ref{max-optimization-problem}) satisfies the KKT conditions
$$
 \left\{\begin{array}{ll}
\nabla \lambda(x)+2\mu x+v=0,\\
v\ge0,\\
v^Tx=0,\\
x\ge0,\\
x^Tx=1.
\end{array}\right.
 $$

Using $v^Tx=0$ and $x^Tg(x)=0$, by taking the dot product with $x$ in the first equation, we get that $\mu=0$. So, the first equation could be written as $v=-\nabla \lambda(x)$. Since $v\ge0$ and $\mathcal{B}$ is positive definite, it follows that
$$(\lambda,x)\in \{\mathbb{R}\times\mathbb{R}^n \backslash \{0\}:0\le x \bot (\lambda\mathcal{B}-\mathcal{A})x^{m-1}\ge0\},$$
i.e. any equilibrium solution $x$ of (\ref{max-optimization-problem}) is a solution of the symmetric TEiCP (\ref{TEiCP.problem}).

Furthermore, the global maximum/minmum of $\lambda (x)$ in $\mathbb{S}^{n-1}_+$ is corresponding to the extreme value of Pareto eigenpair of $(\mathcal{A},\mathcal{B})$ \cite{SQ15-JGO,LHQ2015-COAP1} if $\mathcal{B}$ is strictly copositive, i.e., $\mathcal{B}x^m>0$ for any $x\in \mathbb{R}_+^n \backslash \{0\}$. The concept of copositive tensor is introduced by Qi \cite{Qi2013}. A tensor $\mathcal{A}$ is said copositive if $\mathcal{A}x^m\ge0$ for all $x\in \mathbb{R}_+^n \backslash \{0\}$. $\mathcal{A}$ is copositive (strictly copostive) if and only if all of its Pareto H-eigenvalues or Z-eigenvalues are nonnegative (positve, respectively) \cite{SQ15-JGO}.

\begin{proposition}
Let $\mathcal{A},\mathcal{B} \in \mathbb{S}^{[m,n]}$, and $\mathcal{B}$ is copositive.
Let $$\lambda_{TCP}^{max}=\max\{ \lambda: \exists x \in \mathbb{R}^n \backslash \{0\} \ \mbox{such that}\ (\lambda,x)\in \sigma (\mathcal{A}, \mathcal{B})\},$$
and
$$\lambda_{TCP}^{min}=\min\{ \lambda: \exists x \in \mathbb{R}^n \backslash \{0\} \ \mbox{such that}\ (\lambda,x)\in \sigma (\mathcal{A}, \mathcal{B})\}.$$
Then $\lambda_{TCP}^{max}=\max\{\lambda(x): x \in \mathbb{S}^{n-1}_+\}$ and $\lambda_{TCP}^{min}=\min\{\lambda(x): x \in \mathbb{S}^{n-1}_+\}$.
\end{proposition}

If both $\mathcal{A}$ and $\mathcal{B}$ are symmetric and strictly copositive tensors, then we can use logarithmic function as the merit function in
(\ref{max-optimization-problem}). In such a case,  TEiCP (\ref{TEiCP.problem}) could be reformulated to the following nonlinear optimization problem:

\begin{equation}\label{log-max-optimization-problem}
\max f(x)=\ln(\mathcal{A}x^m)-\ln(\mathcal{B}x^m) \ \ \mbox{subject to} \ \ x \in \mathbb{S}^{n-1}_+.
\end{equation}
Its gradient and Hessian are respectively

\begin{equation}\label{gradientoflog}
g(x)\equiv \nabla f(x)=\frac{m(\mathcal{A}x^{m-1})}{\mathcal{A}x^m}-\frac{m(\mathcal{B}x^{m-1})}{\mathcal{B}x^{m}}.
\end{equation}
and
\begin{align*}\label{Hessianoflog}
H(x)\equiv & \nabla^2 f(x)\\
&=\frac{m(m-1)\mathcal{A}x^{m-2}}{\mathcal{A}x^m}-\frac{m(m-1)\mathcal{B}x^{m-2}}{\mathcal{B}x^m}+
\frac{m^2\mathcal{B}x^{m-1}(\mathcal{B}x^{m-1})^T}{(\mathcal{B}x^m)^2}\\
&-\frac{m^2\mathcal{A}x^{m-1}(\mathcal{A}x^{m-1})^T}{(\mathcal{A}x^m)^2}.
\end{align*}

The Hessian is much simpler than that of Rayleigh quotient function in (\ref{max-optimization-problem}).  If one need to use Hessian for computing  Pareto eigenvalue, the logarithmic merit function may be a favorable choice.

\section{Spectral projected gradient methods}

In this section, the spectral projected gradient (SPG) method is applied to the nonlinear programming problem (\ref{max-optimization-problem}). One main feature of SPG is the spectral choice of step length (also called BB stepsize) along the search direction, originally proposed by Barzilai and Borwein \cite{BBmethod88}.  The Barzilai-Borwein method performs much better than the steepest descent gradient method or projected gradient method in practice \cite{Raydan97,DaiFletcher05,BirginMR00}. Especially, when the objective function is a convex quadratic function and $n=2$, a sequence generated by the BB method converges $R$-superlinearly to the global minimizer \cite{BBmethod88}. For any dimension convex quadratic function, it is still globally convergent \cite{Raydan93} but the convergence is $R$-linear \cite{DaiLiao02}.

We firstly present the following spectral projected gradient method with monotone line search.

\noindent\rule{\textwidth}{1pt}
\underline{\textbf{Algorithm 1: Spectral projected gradient (SPG1) algorithm for TEiCP}} \\
\textbf{Given tensors} $\mathcal{A} \in \mathbb{S}^{[m,n]}$ and $\mathcal{B}\in \mathbb{S}_+^{[m,n]}$, an initial unit iterate $x_0\ge0$, parameter $\rho \in (0,1)$. Let $\epsilon>0$ be the tolerance of termination.   Calculate gradient $g(x_0)$, $\beta_0=1/\|g(x_0)\|$. Set k=0. \\
\textbf{Step 1:} Compute $z_k=P_{\Omega}(x_k+\beta_kg_k)$ and the direction $d_k=z_k-x_k.$ \\
\textbf{Step 2:} If $\|d_k\|=0$ then stop: $\lambda(x)=\frac{\mathcal{A}x^m}{\mathcal{B}x^m}$ is a Pareto eigenvalue, and $x$ is a corresponding Pareto eigenvector of TEiCP. Otherwise, set $\alpha\leftarrow 1$.\\
\textbf{Step 3:} If
\begin{equation}\label{linesearch}
f(x_k+\alpha d_k)\ge f(x_k)+\rho \alpha g_k^Td_k,
\end{equation}
then define $x_{k+1}=x_k+\alpha_k d_k$, $s_k=x_{k+1}-x_k$, $y_k=g_{k+1}-g_k$. Otherwise, set $\alpha=0.5\alpha$ and try again.\\
\textbf{Step 4:} Compute  $b_k=\langle s_k, y_k\rangle $. If $b_k\le0$ set $\beta_{k+1}=\beta_{\max}$; else, compute $a_k=\langle s_k, s_k\rangle$ and
$\beta_{k+1}=\max\{\beta_{\min},\min\{\beta_{\max},\frac{a_k}{b_k}\}\}.$ Set $k :=
k + 1$ and go to \textbf{Step 1}.\\
\noindent\rule{\textwidth}{1pt}

Here $\Omega=\mathbb{S}^{n-1}_+$ is a close convex set. By the projection operation and the convexity of $\Omega$, we know that for all $u\in \Omega$ and $\forall v \in \mathbb{R}^n$,
$$(v-P_{\Omega}(v))^T(u-P_{\Omega}(v))\le 0.$$
Set $u=x$ and $v=x+\beta g(x)$ in the above inequality, then we have
$$\beta g(x)^T[x-P_{\Omega}(x+\beta g(x)]+\|x-P_{\Omega}(x+\beta g(x)\|^2\le0.$$
Let $d_{\beta}(x)=P_{\Omega}(x+\beta g(x))-x$ with $\beta>0$, then we have the following lemma.

\begin{lemma}
 For all $x\in \Omega$, $\beta\in (0,\beta_{\max}]$, we have
\begin{equation}\label{eq:ascentdirection}
g(x)^Td_{\beta}(x)\ge \frac{1}{\beta}\|d_{\beta}(x)\|^2_2\ge \frac{1}{\beta_{\max}}\|d_{\beta}(x)\|^2_2.
\end{equation}
\end{lemma}

From (\ref{eq:ascentdirection}), we know that $d_k$ is an ascent direction. Hence, a stepsize satisfying (\ref{linesearch}) will be found after a finite number of trials, and the SPG algorithm is well defined.
When $\beta=\frac{\langle s_k, s_k\rangle}{\langle s_k, y_k\rangle}$ in $d_{\beta}(x_k)$,  we call it spectral projected gradient (SPG). The vector
$d_{\beta}(x^*)$ vanishes if and only if $x^*$ is a constrained stationary point of optimization problem (\ref{max-optimization-problem})/(\ref{log-max-optimization-problem}).  The convergence of SPG method is established as follows. The proof is similar to that in \cite{BirginMR00}.

\begin{theorem} \label{the:globalconvergence} Let $\{x_{k}\}$ is generated by SPG1 Algorithm. If there is a vector $x_k$ such that $d_{\beta}(x_k)=0$, then $(\lambda (x_k), x_k)$ is a solution of the symmetric TEiCP. Otherwise,  any accumulation point of the sequence $\{x_{k}\}$ is a constrained stationary point, i.e., the sequence $\{\lambda (x_{k})\}$ converges to a Pareto eigenvalue of the symmetric TEiCP.
\end{theorem}

\noindent{\bf Proof.} Let $x^*$ be an accumulation point of $\{x_{k}\}$, and relabel $\{x_{k}\}$ a subsequence converging to $x^*$. According to the Proposition 2, we just need to show that $x^*$ is a constrained stationary point of the optimization problem. Let us suppose by way of contradiction that $x^*$ is not a constrained stationary point. So, by continuity and compactness, there exist $\delta>0$ such that $\|d_{\beta}(x^*)\|\ge\delta>0 $ for all $\beta\in (0,\beta_{\max}]$. Furthermore, using the Lemma 1, we have $g(x^*)^Td_{\beta}(x^*)\ge \frac{\delta^2}{\beta_{\max}}$ for all $\beta\in (0,\beta_{\max}]$, which implies that for $k$ larger enough on the subsequence that converges to $x^*$, $g(x_k)^Td_{\beta}(x_k)>c$ for all $\beta\in [\beta_{\min},\beta_{\max}]$. Here, we can set $c=\frac{\delta^2}{2\beta_{\max}}>0$. We consider two cases.

Firstly, assume that $\inf \alpha_k\ge \varepsilon >0$. By continuity, for sufficiently large $k$, $\|d_{\beta_k}(x_k)\|\ge\delta/2$. From the line search condition (\ref{linesearch}), we have
$$f(x_k+\alpha d_k)-f(x_k)\ge \rho \alpha g_k^Td_k\ge \frac{\rho \varepsilon \delta^2}{4 \beta_{\max}}. $$
Clearly,  when $k\rightarrow \infty$, $f(x_k)\rightarrow \infty$, which is a contradiction. In fact, $f$ is a continuous function and so $f(x_k)\rightarrow f(x^*)$.

Assume that $\inf \alpha_k=0$.
Since $\inf \alpha_k=0$, there exists a subsequence $\{x_{k}\}_K$ such that $\lim_{k\in K}\alpha_k=0$. In such a case, from the way $\alpha_k$ is chosen in (\ref{linesearch}), there exists an index $\bar{k}$ sufficiently large such that for all $k\ge \bar{k}$, $k\in K$, for which $\alpha_k/0.5$ fails to satisfy condition (\ref{linesearch}), i.e.,
$f(x_k+2\alpha d_k)-f(x_k)<2 \rho \alpha_k g_k^Td_k$. Hence,
$$\frac{f(x_k+2\alpha d_k)-f(x_k)}{2 \alpha_k}<\rho g_k^Td_k.$$
By the mean value theorem, we can rewrite this relation as
$$d_k^Tg(x_k+t_kd_k)<\rho g_k^Td_k \; \mbox{for all} \; k\in K, \; k\ge \bar{k},$$
where $t_k\in [0, 2\alpha_k]$ that goes to zero as $k\in K$ goes to infinity. Taking limits in the above inequality, we deduce that $(1-\rho)g(x^*)^Td(x^*)\le0$. Since $1-\rho>0$ and $g_k^Td_k>0$ for all k, then $g(x^*)^Td(x^*)=0.$ By continuity, this indicates that for $k$ large enough on the subsequence we have that $g_k^Td_k<c/2$, which contradicts to $g_k^Td_k>c$.

Therefore,  any accumulation point of the sequence $\{x_{k}\}$ is a constrained stationary point. By using the Proposition 2, it follows that the sequence $\{\lambda (x_{k})\}$ converges to a Pareto eigenvalue of the symmetric TEiCP.  $\quad\Box$\\

In the rest of this section, we would like to present the following SPG algorithm for TEiCP with curvilinear search. Its global convergence could be established similarly.

\noindent\rule{\textwidth}{1pt}
\underline{\textbf{Algorithm 2: Spectral projected gradient (SPG2) algorithm for TEiCP}} \\
\textbf{Given tensors} $\mathcal{A} \in \mathbb{S}^{[m,n]}$ and $\mathcal{B}\in \mathbb{S}_+^{[m,n]}$, an initial unit iterate $x_0\ge0$, parameter $\rho \in (0,1)$. Let $\epsilon>0$ be the tolerance of termination. Calculate gradient $g(x_0)$, $\beta_0=1/\|g(x_0)\|$. Set k=0. \\
\textbf{Step 1:} If  $\|P_{\Omega}(x_k+\beta_kg_k)-x_k\|<\epsilon$, stop, declaring $\lambda(x)=\frac{\mathcal{A}x^m}{\mathcal{B}x^m}$ is a Pareto eigenvalue, and $x$ is a corresponding Pareto eigenvector of TEiCP.\\
\textbf{Step 2:} Set $\alpha \leftarrow \beta_k$.  \\
\textbf{Step 3:} Set $x_+=P_{\Omega}(x_k+\alpha g_k).$ If
\begin{equation}\label{linesearch2}
f(x_+)\ge f(x_k)+\rho \alpha g_k^T(x_+-x_k),
\end{equation}
then define $x_{k+1}=x_+$, $s_k=x_{k+1}-x_k$, $y_k=g_{k+1}-g_k$. Otherwise, set $\alpha=0.5\alpha$ and try again.\\
\textbf{Step 4:} Compute  $b_k=\langle s_k, y_k\rangle $. If $b_k\le0$ set $\beta_{k+1}=\beta_{\max}$; else, compute $a_k=\langle s_k, s_k\rangle$ and
$$\beta_{k+1}=\max\{\beta_{\min},\min\{\beta_{\max},\frac{a_k}{b_k}\}\}.$$ Set $k :=
k + 1$ and go to \textbf{Step 1}.\\
\noindent\rule{\textwidth}{1pt}



\section{Numerical experiments}

In this section, we present some numerical results to illustrate the effectiveness of the spectral projected gradient (SPG) methods, which were compared with the Scaling-and-Projection Algorithm (SPA) proposed by Ling, He and Qi \cite{LHQ2015-COAP1} and the shifted projected power (SPP) method for TEiCP proposed in \cite{ChenYangYe15}.

Both SPG1 and SPG2 are monotone ascent method. $\{x_k\}$ are always located in the feasible region $\Omega$. In general, the merit function $f(x)$ is chosen to be the Rayleigh quotient function in (\ref{max-optimization-problem}). In the line search procedure of the SPG1 method, we used the one-dimensional quadratic interpolation to compute the stepsize $\alpha$ such as
$$\alpha\leftarrow\frac{-\alpha^2 g_k^Td_k} { 2(f(x_k+\alpha d_k) - f(x_k) - \alpha g_k^Td_k)}.$$
In the implementation, we terminate the algorithm once
$$\|x_{k+1}-x_k\|\le \epsilon,\; \mbox{or} \;\|g(x_k)\|\le \epsilon,\; \mbox{or} \; \|\lambda_{k+1}-\lambda_k\|\le \epsilon. $$
We accept $\epsilon=10^{-6},$ and set the parameter $\rho=10^{-4}$, $\beta_{\max}=\frac{1}{\|g_k\|}$ and $\beta_{\min}=\|g_k\|$. For SPP and SSPA, the parameter $\tau=0.05$. In all numerical experiments, the maximum iterations is 500.
The experiments were done on a laptop with Intel Core 2 Duo CPU with a 4GB RAM, using MATLAB R2014b, and the Tensor Toolbox \cite{TensorToolbox}.

We firstly describe the so-called shifted projected power (SPP) algorithm and the scaling-and-projection algorithm (SPA) as follows.

\noindent\rule{\textwidth}{1pt}
\underline{\textbf{Algorithm 3: Shifted Projected Power (SPP) algorithm \cite{ChenYangYe15} }} \\
\textbf{Given tensors} $\mathcal{A} \in \mathbb{S}^{[m,n]}$ and $\mathcal{B}\in \mathbb{S}_+^{[m,n]}$, an initial unit iterate $x_0\ge0$.
 Let $\epsilon>0$ be the tolerance on termination. Let $\tau>0$ be the tolerance on being positive definite.\\
\textbf{for} $k=0,1,\ldots$ \textbf{do}\\
\textbf{1:} Compute the gradient $g(x_k)=\nabla f(x_k)$ and the Hessian $H(x_k)=\nabla^2 f(x_k)$, respectively. Let $r_k\leftarrow  \max\{0,(\tau-\lambda_{min}(H_k))/m$, $\nabla \hat{f}(x_k)=\nabla f(x_k)+r_k m x_k$.\\
\textbf{2:} Let $\nabla \hat{f}_+(x_k)=\left\{\begin{array}{ll}
0,& \ \mbox{if} \ (\nabla \hat{f}(x_k))_i<0 ,\\
(\nabla \hat{f}(x_k))_i,& \mbox{otherwise},\\
\end{array}\right.$\\
\textbf{3:} If $\|\nabla \hat{f}_+(x_k)\|\le \epsilon$, stop. Otherwise, $x_{k+1}\leftarrow \nabla \hat{f}_+(x_k)/\|\nabla \hat{f}_+(x_k)\|$. Set k=k+1 and go back to Step 1.\\
\textbf{End for}\\
\noindent\rule{\textwidth}{1pt}

\noindent\rule{\textwidth}{1pt}
\underline{\textbf{Algorithm 4: Scaling-and-Projection Algorithm (SPA) \cite{LHQ2015-COAP1} }} \\
\textbf{Given tensors} $\mathcal{A} \in \mathbb{S}^{[m,n]}$ and $\mathcal{B}\in \mathbb{S}_+^{[m,n]}$. For an initial point  $u_0\ge0$, define
$x_0=u_0/\sqrt[m]{\mathcal{B}(u_0)^m}$. Let $\epsilon>0$ be the tolerance on termination. \\
\textbf{for} $k=0,1,\ldots$ \textbf{do}\\
\textbf{1:} Compute $\lambda_k=\frac{\mathcal{A}(x_k)^m}{\mathcal{B}(x_k)^m}$, the gradient $g(x_k)=\nabla f(x_k)=\mathcal{A}(x_k)^{m-1}-\lambda_k \mathcal{B}(x_k)^{m-1}$. \\
\textbf{2:} If $\|g(x_k)\|\le \epsilon$, stop. Otherwise, let $\alpha_k=\|g_k\|$, compute $u_{k}=P_{\Omega}(x_k+\alpha_kg_k)$, and
 $x_{k+1}=u_k/\sqrt[m]{\mathcal{B}(u_k)^m}$.
 Set k=k+1 and go back to Step 1.\\
\textbf{End for}\\
\noindent\rule{\textwidth}{1pt}

Since the stepsize $\alpha_k$ in SPA approaches to zero as the sequence $\{x_k\}$ gets close to a solution of TEiCP, as shown in \cite{LHQ2015-COAP1}, the number of iterations will increase significantly. In order to improve the efficiency of SPA method, they try to amplify the stepsize and proposed a modification of SPA such as $u_{k}=P_{\Omega}(x_k+s\alpha_kg_k)$ with $s\in (1,8)$ being a constant parameter.  A suitable choice $s$ will get an improvement. But, how to choose it? Anyway, the stepsize $s\alpha_k$ also approaches to zero when the sequence $\{x_k\}$ gets close to a solution of TEiCP. When the merit function $f(x)$ is (locally) convex, this situation will be better. So, we present the following shifted SPA method, in which an adaptive shift could force the objective to be (locally) convex \cite{KoldaMayo14}.

\noindent\rule{\textwidth}{1pt}
\underline{\textbf{Algorithm 5: Shifted Scaling-and-Projection Algorithm (SSPA) }} \\
\textbf{Given tensors} $\mathcal{A} \in \mathbb{S}^{[m,n]}$ and $\mathcal{B}\in \mathbb{S}_+^{[m,n]}$. For an initial point  $u_0\ge0$, define
$x_0=u_0/\sqrt[m]{\mathcal{B}(u_0)^m}$. Compute $\lambda (x_0)=\frac{\mathcal{A}(x_0)^m}{\mathcal{B}(x_0)^m}$. Let $\epsilon>0$ be the tolerance on termination. Let $\tau>0$ be the tolerance on being positive definite.\\
\textbf{for} $k=0,1,\ldots$ \textbf{do}\\
\textbf{1:}  Compute $y(x_k)=\mathcal{A}(x_k)^{m-1}-\lambda_k \mathcal{B}(x_k)^{m-1}$, the Hessian $H(x_k)=\nabla^2 f(x_k)$, respectively. Let $r_k\leftarrow  \max\{0,(\tau-\lambda_{min}(H_k))/m$, $\hat{g}(x_k)=\nabla \hat{f}(x_k)=y(x_k)+r_k m x_k$.\\
\textbf{2:} Let $\alpha_k=\|\hat{g}(x_k)\|$, compute $u_{k}=P_{\Omega}(x_k+\alpha_k\hat{g}_k)$, and
 $x_{k+1}=u_k/\sqrt[m]{\mathcal{B}(u_k)^m}$, $\lambda (x_{k+1})=\frac{\mathcal{A}(x_{k+1})^m}{\mathcal{B}(x_{k+1})^m}$.\\
\textbf{3:} If $|\lambda (x_{k+1})-\lambda (x_k)|\le \epsilon$, stop. Otherwise,
 Set k=k+1 and go back to Step 1.\\
\textbf{End for}\\
\noindent\rule{\textwidth}{1pt}

\subsection{Comparison with SPA for computing Pareto Z-eigenpairs}
The following example is originally from \cite{KR02} and was used in evaluating the SS-HOPM algorithm in \cite{KoldaMayo11} and the GEAP algorithm in \cite{KoldaMayo14} for computing Z-eigenpairs.

{\it Example 1} (Kofidis and Regalia \cite{KR02}). Let $\mathcal{A} \in \mathbb{S}^{[4,3]}$ be the symmetric tensor defined by
\begin{equation*}
\begin{split}
a_{1111}=0.2883, \ \ a_{1112}=-0.0031,\ \  a_{1113}=0.1973, \ \ a_{1122}=-0.2485,&\\
a_{1223}= 0.1862, \ \ \  a_{1133}=\ 0.3847,\ \  a_{1222}=0.2972, \ \ a_{1123}=-0.2939,&\\
a_{1233}=0.0919, \ \ a_{1333}=-0.3619,\ \  a_{2222}=0.1241, \ \ a_{2223}= -0.3420,&\\
a_{2233}=0.2127, \ \ a_{2333}=0.2727,\ \  a_{3333}=-0.3054.&
\end{split}
\end{equation*}

\begin{figure}\label{fig:SPG-SPA-Z-eigen-1}
$$
\begin{array}{cc}
\centerline{\includegraphics[width=1.1\textwidth]{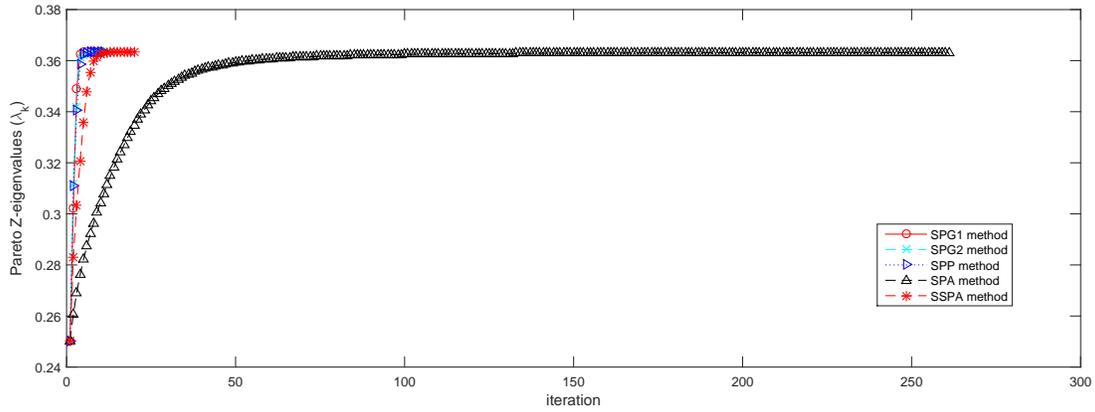}}
\end{array}
$$
 \caption{Comparison with SPA algorithm for computing Pareto Z-eigenvalues of $\mathcal{A}$ from {\bf \it Example 1}, and the starting point is $x_0=[1.0;1.0;1.0]$ }

\end{figure}

\begin{figure}\label{fig:SPG-SPP-SSPA-Z-eigen-100-1}
$$
\begin{array}{cc}
\centerline{\includegraphics[width=1.1\textwidth]{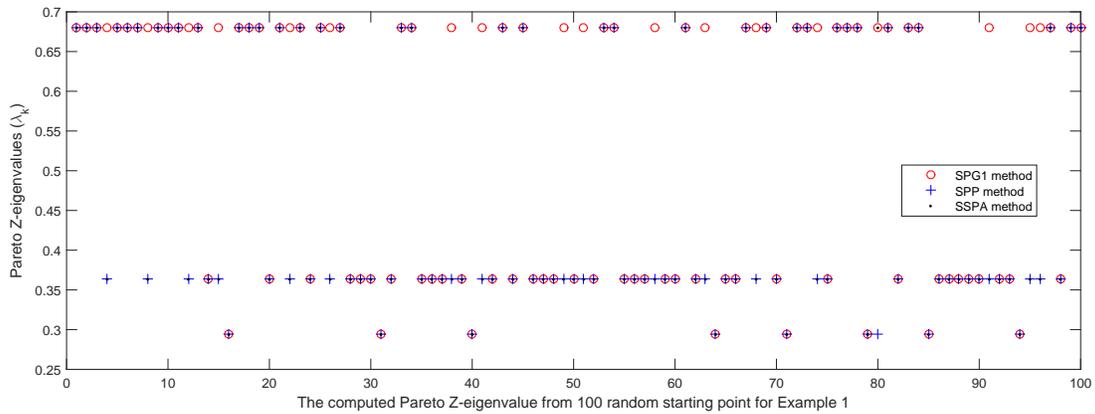}}
\end{array}
$$
 \caption{The computed Pareto Z-eigenvalues by SPG, SPP and SSPA in the 100 runs on the $\mathcal{A}$ from {\bf \it Example 1}. }

\end{figure}

To compare the convergence in terms of the number of iterations. Figure 1 shows the results for computing Pareto Z-eigenvalues of $\mathcal{A}$ from {\bf \it Example 1}, and the starting point is $x_0=[1.0;1.0;1.0]$. In this case, all of the SPG1, SPG2, SPP, SSPA can reach the same Pareto Z-eigenvalue 0.3633. SPG1 method just need run 9 iterations in 0.1716 seconds while SPA method need run 260 iterations in  3.3696 seconds. SPP is similar to SPG1 method in this case. SPG2 need run 13 iterations in 0.4368 seconds. As we can see, comparing with SPA method, SSPA method get a great improvement. SSPA method just need run 19 iterations in 0.2964 seconds.

\begin{center}
\small{Table 1. Comparison results for computing Pareto Z-eigenvalues of $\mathcal{A}$ from {\bf \it Example 1}.}\\
\begin{tabular}{|c|c|c|c|c|c|}
\hline
Alg. & $\lambda$ & Eigenvector &  Its. & Error &  Time (sec.)\\
\hline
SPG1& 0.3633 &  [0.2678;0.6446;0.7161] & 9 &  5.43e-07 & 0.1716 \\
SPG2& 0.3633 &  [0.2677;0.6445;0.7162] & 13 &  3.94e-08 &  0.4368 \\
SPP& 0.3633 & [0.2679;0.6448;0.7158] & 10 &  5.46e-07 & 0.1404 \\
SPA& 0.3632 & [0.2771;0.6461;0.7112] & 260 &  9.96e-07 & 3.3696 \\
SSPA& 0.3633 & [0.2683;0.6449;0.7156] & 19 &  9.00e-07 & 0.2964 \\
\hline
\end{tabular}
\end{center}

{\it Example 2}.  Let $\mathcal{A} \in \mathbb{S}^{[4,n]}$ be the diagonal tensor defined by
$a_{iiii}=\frac{i-1}{i},$ for $i=1,\ldots,n$.

\begin{center}
\small{Table 2. Comparison results for computing Pareto Z-eigenvalues of $\mathcal{A}$ from {\bf \it Example 2} with $n=5$.}\\
\begin{tabular}{|c|c|c|c|c|c|}
\hline
Alg. & $\lambda$ & Eigenvector &  Its. & Error &  Time (sec.)\\
\hline
SPG1& 0.8 &  [0;0;0;0;1] & 3 &  0.0* & 0.0312 \\
SPG2& 0.8 &  [0;0;0;0;1] & 4 &  0.0* &  0.0312 \\
SPP&  0.8 & [0;2.95e-10;7.13e-09;3.22e-07;0.9999] & 7 &  7.02e-10 & 0.0624 \\
SPA&  0.7999 & [0.0014;0.0024;0.0037;0.0063;0.9999] & 286 &  9.87e-07 & 3.4008 \\
SSPA& 0.8 & [8.86e-05;1.51e-04;2.38e-04;4.31e-04;0.9999] & 60 &  7.45e-07 & 0.7800 \\
\hline
\end{tabular}
\end{center}

Figure 3 shows the results for computing Pareto Z-eigenvalues of $\mathcal{A}$ from {\bf \it Example 2} with $n=5$, and the starting point is $x_0=[1.0;1.0;1.0;1.0;1.0]$. In this case, all of the SPG1, SPG2, SPP, SPA, SSPA can reach the largest Pareto Z-eigenvalue 0.8. SPG1 just need 3 iterations in 0.0312 seconds while SPA need run 286 iterations in 3.4008 seconds.  Comparing with SPA method, SSPA method get a great improvement again. But, SSPA method is still slower than the other three methods in this case.

{\it Example 3.} Let $\mathcal{A} \in \mathbb{S}^{[4,3]}$ be the symmetric tensor defined by: Firstly, set $A=tensor(zeros(3,3,3,3)),$ and
\begin{equation*}
\begin{split}
a_{1111}=1.00397, \ \ a_{2222}=0.99397,\ \  a_{3333}=1.00207,&\\
a_{1222}= 0.00401, \ \ \  a_{2111}=\ 0.00788,\ \  a_{3111}=0.00001,&\\
a_{3222}=0.00005, \ \ a_{1333}=0.99603,\ \  a_{2333}=1.0040,&\\
\end{split}
\end{equation*}
and then using $A = symmetrize(A)$ to symmetrize it.

\begin{figure}\label{fig:SPG-SPA-Z-eigen-2}
$$
\begin{array}{cc}
\centerline{\includegraphics[width=1.1\textwidth]{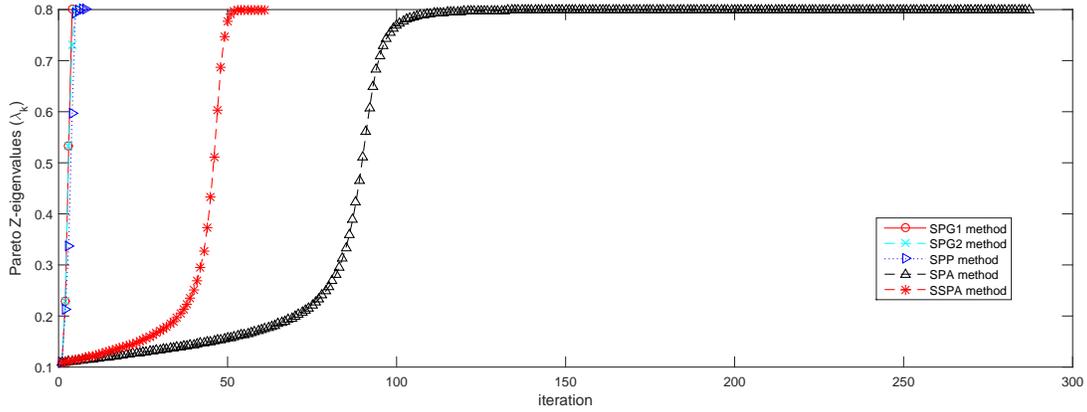}}
\end{array}
$$
 \caption{Comparison with SPA algorithm for computing Pareto Z-eigenvalues of $\mathcal{A}$ from {\bf \it Example 2} with $n=5$  and the starting point is $x_0=[1.0;1.0;1.0;1.0;1.0]$ }

\end{figure}

\begin{figure}\label{fig:SPG-SPP-SSPA-Z-eigen-100-2}
$$
\begin{array}{cc}
\centerline{\includegraphics[width=1.1\textwidth]{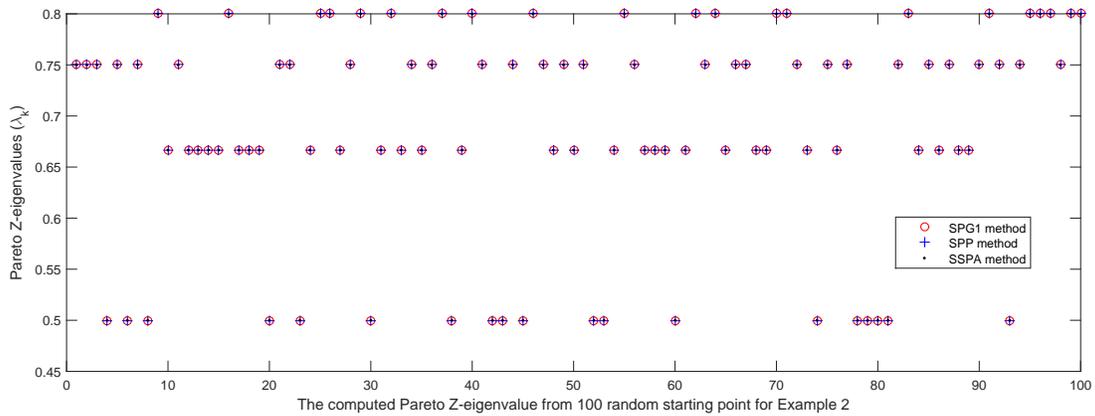}}
\end{array}
$$
 \caption{The computed Pareto Z-eigenvalues by SPG, SPP and SSPA in the 100 runs on the $\mathcal{A}$ from {\bf \it Example 2} with $n=5$. }

\end{figure}

\begin{figure}\label{fig:SPG-SPA-Z-eigen-3}
$$
\begin{array}{cc}
\centerline{\includegraphics[width=1.1\textwidth]{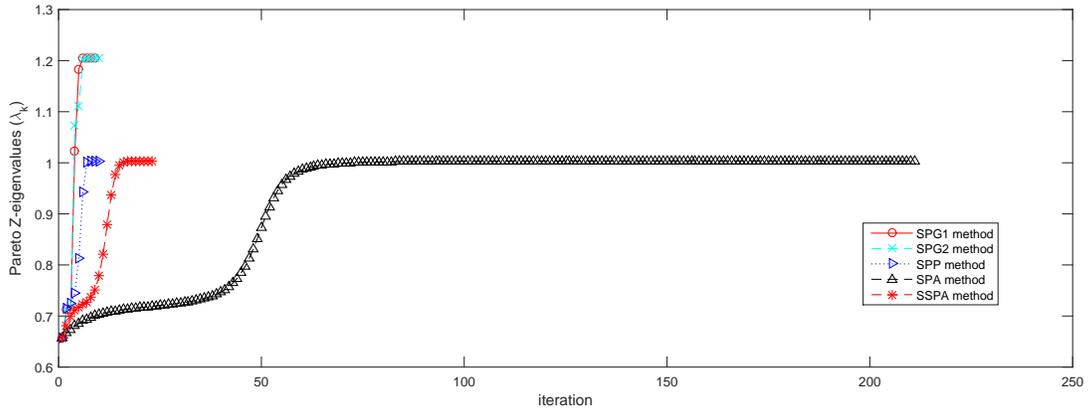}}
\end{array}
$$
 \caption{Comparison with SPA algorithm for computing Pareto Z-eigenvalues of $\mathcal{A}$ from {\bf \it Example 3}, and the starting point is $x_0=[0.9015;0.3183;0.5970]$ }

\end{figure}

\begin{figure}\label{fig:SPG-SPP-SSPA-Z-eigen-100-3}
$$
\begin{array}{cc}
\centerline{\includegraphics[width=1.1\textwidth]{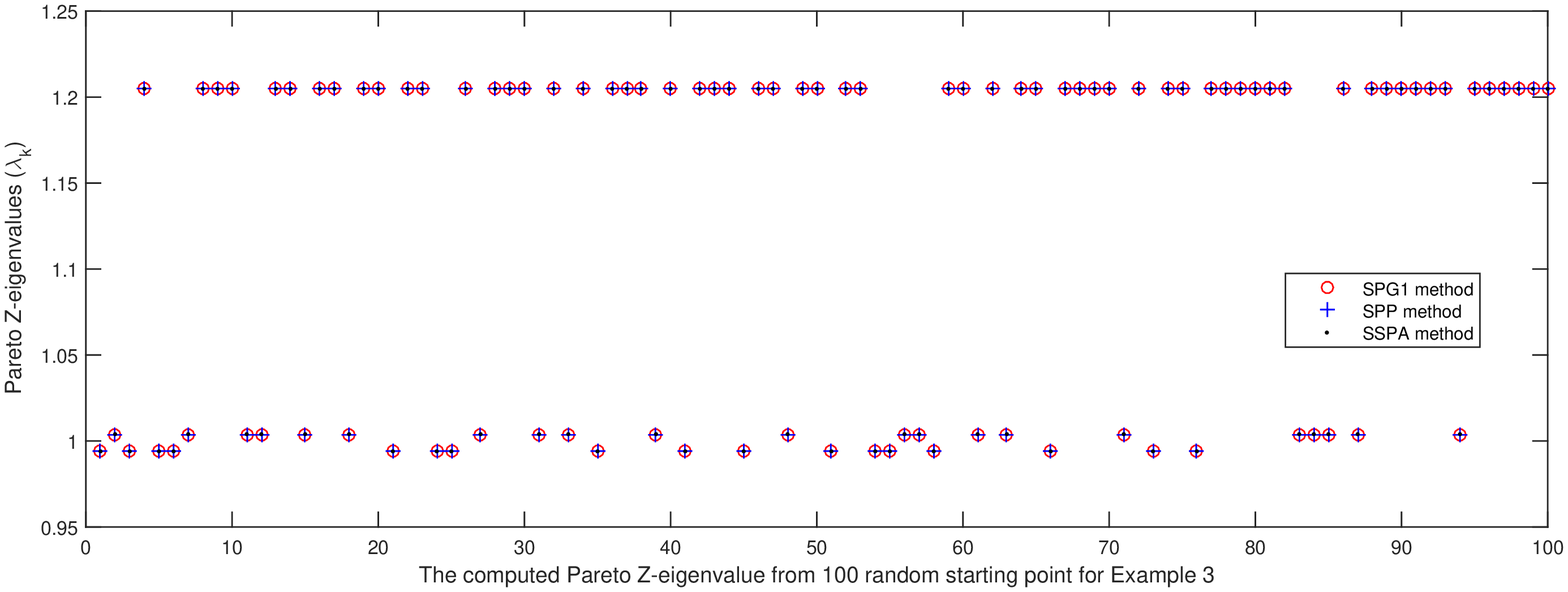}}
\end{array}
$$
 \caption{The computed Pareto Z-eigenvalues by SPG, SPP and SSPA in the 100 runs on the $\mathcal{A}$ from {\bf \it Example 3}. }

\end{figure}

\begin{center}
\small{Table 3. Comparison results for computing Pareto Z-eigenvalues of $\mathcal{A}$ from {\bf \it Example 3}.}\\
\begin{tabular}{|c|c|c|c|c|c|}
\hline
Alg. & $\lambda$ & Eigenvector &  Its. & Error &  Time (sec.)\\
\hline
SPG1& 1.2048 &  [0.1905;0.1920;0.9627] & 8 &  5.25e-07 & 0.1248 \\
SPG2& 1.2048 &  [0.1902;0.1918;0.9628] & 9 &  1.71e-07 &  0.2028 \\
SPP& 1.0040 & [1.0;0.0020;2.72e-06] & 9 &  7.37e-10 & 0.1560 \\
SPA& 1.0039 & [1.0;0.0026;0.0062] & 210 &  9.81e-07 & 2.6832 \\
SSPA& 1.0040 & [1.0;0.0020;2.94e-04] & 22 &  5.19e-07 & 0.2652 \\
\hline
\end{tabular}
\end{center}

Figure 3 shows the results for computing Pareto Z-eigenvalues of $\mathcal{A}$ from {\bf \it Example 3}, and the starting point is $x_0=[0.9015;0.3183;0.5970]$. In this case, both SPG1 and SPG2 reach the largest Pareto Z-eigenvalue $\lambda=1.2048$ while the other three methods reach the Z-eigenvalue $\lambda=1.0040$.  SPG and SPP need 8/9 iterations while SPA need run 210 iterations, and SSPA need 22 iterations in this case.

We also used 100 random starting guesses, each entry selected uniformly at random from the interval $[0,1]$, to test SPG1,SPP and SSPA, respectively.  For each set of experiments, the same set of random starts was used. We listed the median number of iterations until convergence, and the average run time in the 100 experiments in Table 4. The computed Pareto Z-eigenvalues were listed in Fig.2, Fig.4, Fig.6 for Example 1,2,3, respectively. As we can see, most of time, all of the three methods can reach the same Pareto Z-eigenvalue. But for Example 1, it is seems that SPG1 method could reach the largest Pareto Z-eigenvalue with a higher probability.

\begin{center}
\small{Table 4. Comparison results for 100 random test on computing Pareto Z-eigenvalues of $\mathcal{A}$ from {\bf \it Ex.1,Ex.2($n=5$) and Ex.3}.}\\
\begin{tabular}{|c|c|c|c|c|c|c|c|c|}
\hline \multicolumn{ 1}{|c|}{Algorithm}  &  \multicolumn{ 2}{|c|}{SPG1} &  \multicolumn{ 2}{|c|}{SPP} &  \multicolumn{ 2}{|c|}{SSPA} \\
\cline{2-7}
\multicolumn{ 1}{|c|}{Example} & Its. &  Time & Its. &  Time & Its. &  Time\\
\hline
Ex. 1 &  7.41 & 0.1407 &  8.17 & 0.1259&  15.61 & 0.2253 \\
Ex. 2 &  2.11 & 0.0381 &  5.21 & 0.0841&  37.08 & 0.5129 \\
Ex. 3 &  4.79 & 0.0894 &  5.20 & 0.0853&  14.30 & 0.2058 \\
\hline
\end{tabular}
\end{center}

\subsection{Comparison with SPP for computing Pareto H-eigenpairs}
In this subsection, we test SPG1, SPG2, and SPP method for finding Pareto H-eigenpairs of $\mathcal{A}$ from Examples 4-6 ($n=5$):

{\it Example 4} (Nie and Wang \cite{NW14}). Let $\mathcal{A} \in \mathbb{S}^{[4,n]}$ be the symmetric tensor defined by
$$a_{ijkl}=\sin(i+j+k+l) \ \ (1\le i,j,k,l\le n).$$

{\it Example 5} (Nie and Wang \cite{NW14}). Let $\mathcal{A} \in \mathbb{S}^{[4,n]}$ be the symmetric tensor defined by
$$a_{ijkl}=\tan(i)+\tan(j)+\tan(k)+\tan(l) \ \ (1\le i,j,k,l\le n).$$

{\it Example 6} (Nie and Wang \cite{NW14}).  Let $\mathcal{A} \in \mathbb{S}^{[4,n]}$ be the tensor defined by
$$a_{ijkl}=\frac{(-1)^i}{i}+\frac{(-1)^j}{j}+\frac{(-1)^k}{k}+\frac{(-1)^l}{l}, (1\le i,j,k,l \le n). $$

\begin{figure}\label{fig:SPG-SPP-H-eigen-4}
$$
\begin{array}{cc}
\centerline{\includegraphics[width=1.1\textwidth]{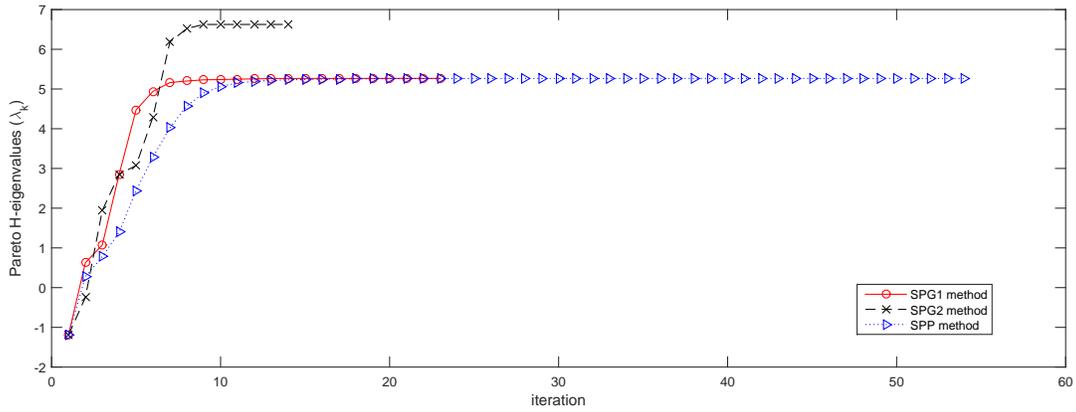}}
\end{array}
$$
 \caption{Comparison SPG with SPP algorithm for computing H-eigenvalues of $\mathcal{A}$ from {\bf \it Example 4} (n=5), and the starting point is $x_0=[0.3319;0.8397;0.3717;0.8282;0.1765]$. }

\end{figure}

\begin{figure}\label{fig:SPG-SPP-H-eigen-100-4}
$$
\begin{array}{cc}
\centerline{\includegraphics[width=1.1\textwidth]{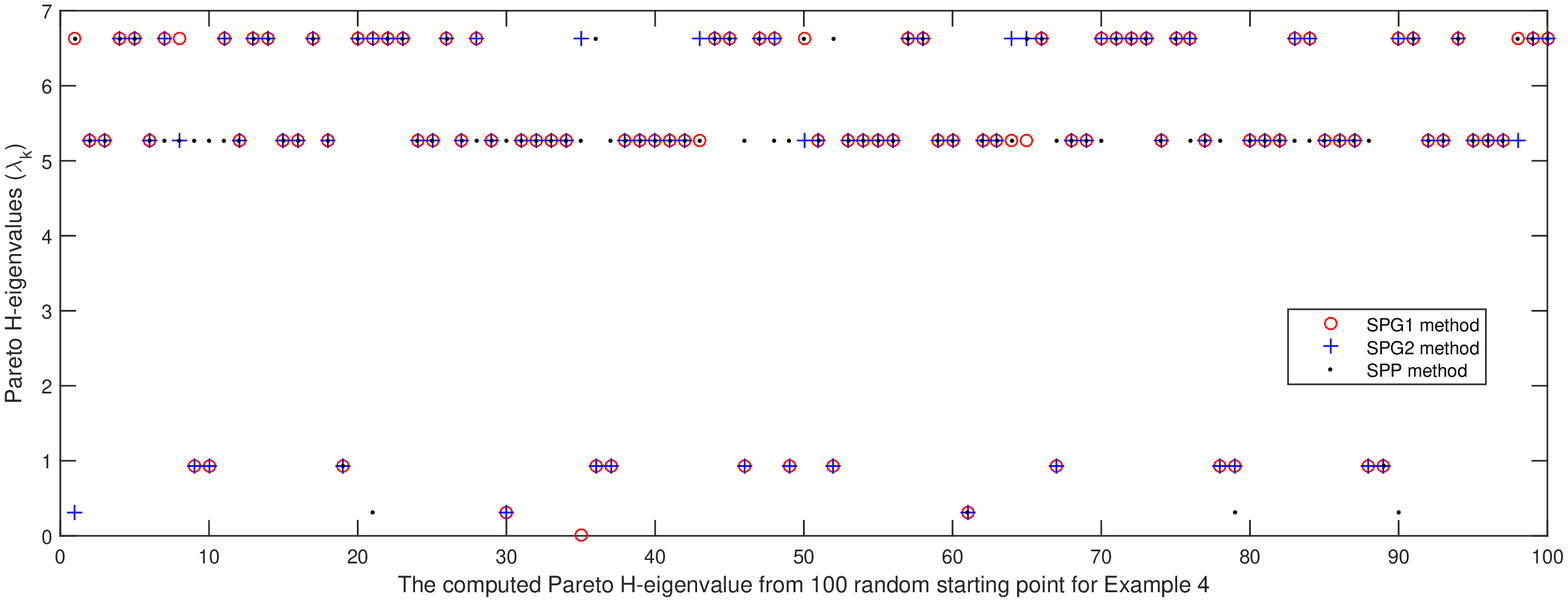}}
\end{array}
$$
 \caption{The computed Pareto H-eigenvalues by SPG1, SPG2, and SPP in the 100 runs on the $\mathcal{A}$ from {\bf \it Example 4}. }

\end{figure}

\begin{figure}\label{fig:SPG-SPP-H-eigen-5}
$$
\begin{array}{cc}
\centerline{\includegraphics[width=1.1\textwidth]{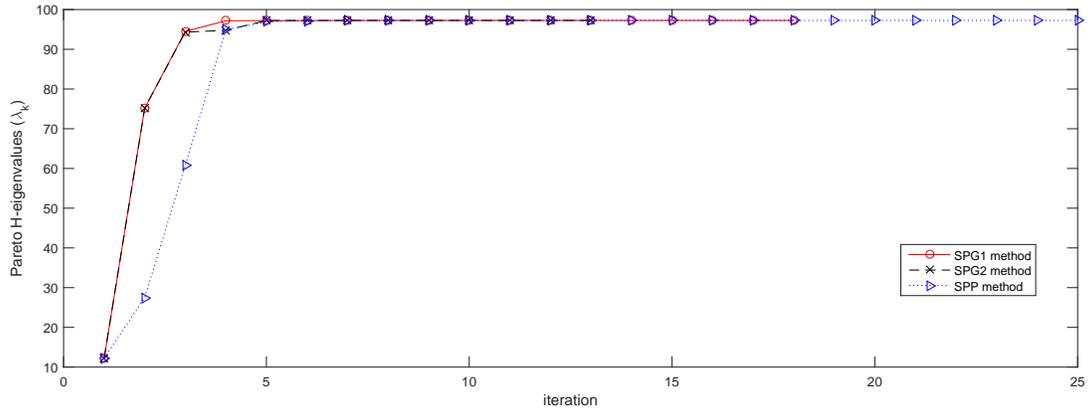}}
\end{array}
$$
 \caption{Comparison SPG with SPP algorithm for computing H-eigenvalues of $\mathcal{A}$ from {\bf \it Example 5} (n=5), and the starting point is $x_0=[0.2291;0.0922;0.2409;0.9025;0.21734]$. }

\end{figure}

\begin{figure}\label{fig:SPG-SPP-H-eigen-100-5}
$$
\begin{array}{cc}
\centerline{\includegraphics[width=1.1\textwidth]{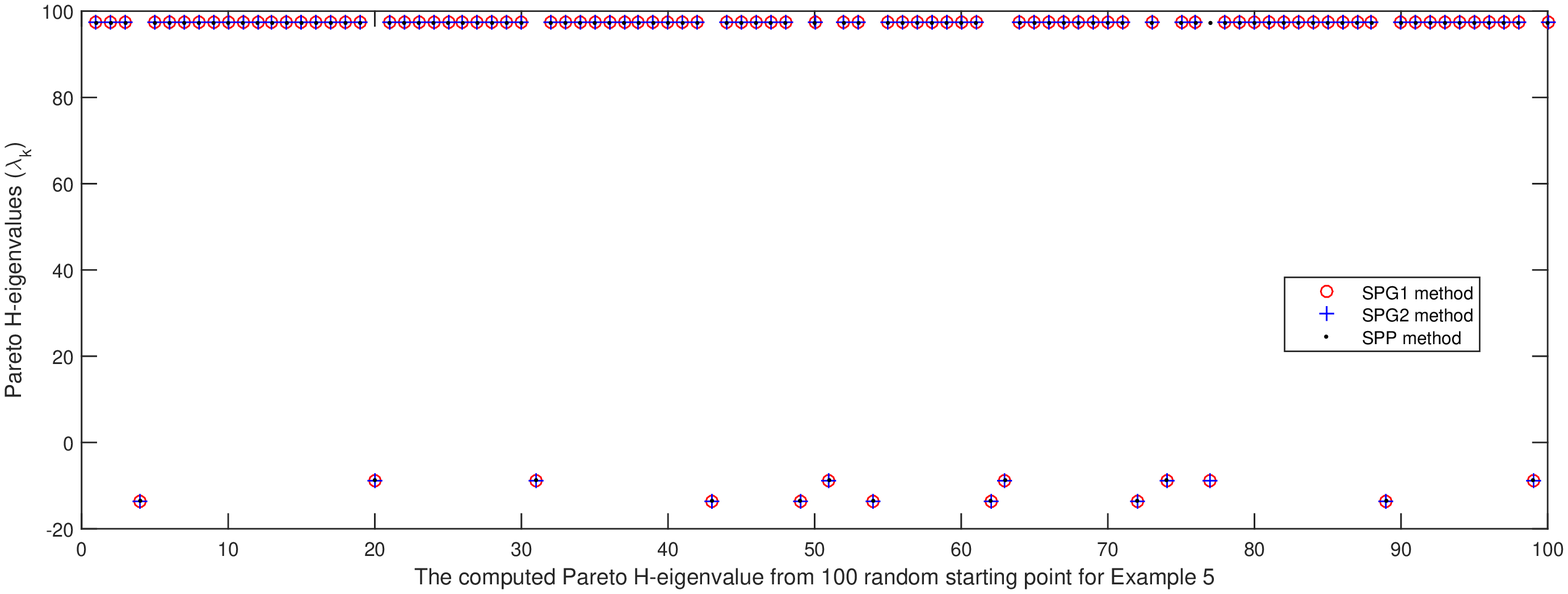}}
\end{array}
$$
 \caption{The computed Pareto H-eigenvalues by SPG1, SPG2, and SPP in the 100 runs on the $\mathcal{A}$ from {\bf \it Example 5}. }

\end{figure}

\begin{figure}\label{fig:SPG-SPP-H-eigen-6}
$$
\begin{array}{cc}
\centerline{\includegraphics[width=1.1\textwidth]{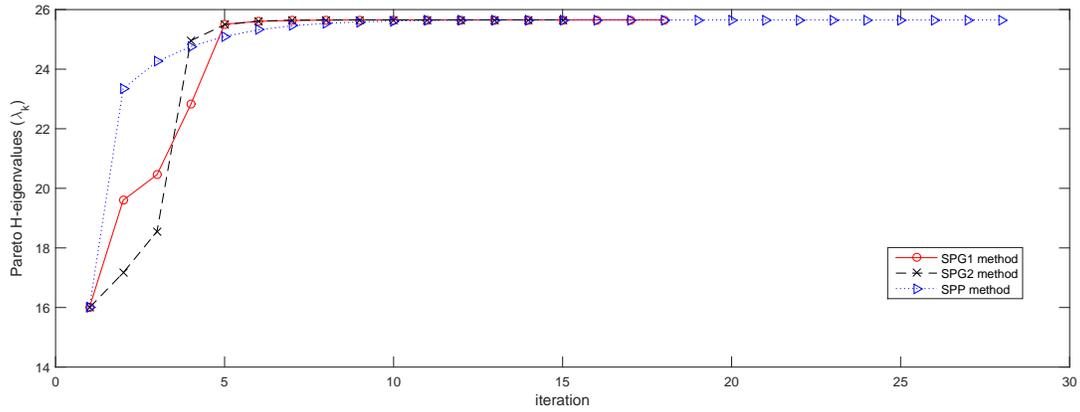}}
\end{array}
$$
 \caption{Comparison SPG with SPP algorithm for computing H-eigenvalues of $\mathcal{A}$ from {\bf \it Example 6} (n=5), and the starting point is $x_0=[0.1846;0.8337;0.1696;0.9532;0.7225]$. }

\end{figure}

\begin{figure}\label{fig:SPG-SPP-H-eigen-100-6}
$$
\begin{array}{cc}
\centerline{\includegraphics[width=1.1\textwidth]{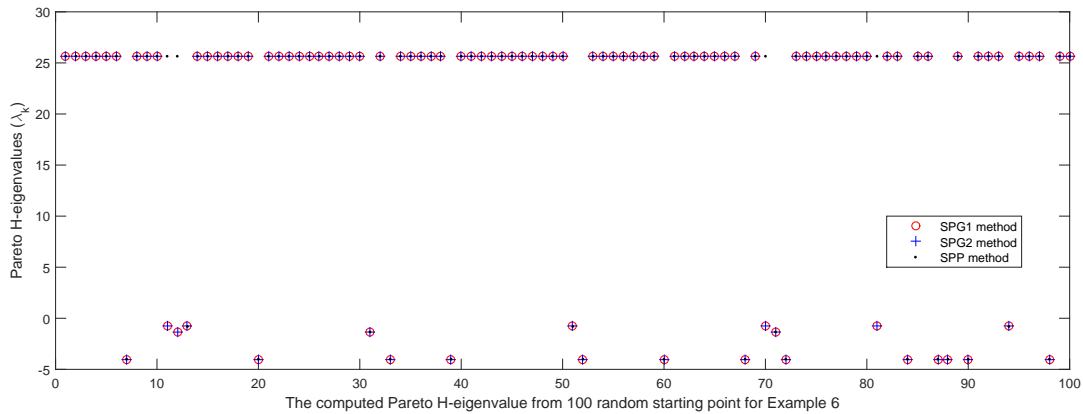}}
\end{array}
$$
 \caption{The computed Pareto H-eigenvalues by SPG1, SPG2, and SPP in the 100 runs on the $\mathcal{A}$ from {\bf \it Example 6}. }

\end{figure}

To compare the convergence in terms of the number of iterations. Figure 7 shows the results for computing Pareto H-eigenvalues of $\mathcal{A}$ from {\bf \it Example 4}, and the starting point is $x_0=[0.3319;0.8397;0.3717;0.8282;0.1765]$. In this case, both of SPG1 method and SPP method can find the same Pareto H-eigenvalue 5.2664 while SPG2 method finding the Pareto H-eigenvalue 6.6255. SPG1 method need run 22 iterations in 0.4368 seconds while SPP method need run 53 iterations in  0.6240 seconds. SPG2 method just need run 13 iterations in 0.4524 seconds.

Figure 9 shows the results for computing Pareto H-eigenvalues of $\mathcal{A}$ from {\bf \it Example 5}, and the starting point is $x_0=[0.2291;0.0922;0.2409;0.9025;0.21734]$. In this case, all of the three methods can find the largest Pareto H-eigenvalue 97.2637. SPG1 method need run 17 iterations in 0.3900 seconds while SPG2 method need run 12 iterations in  0.3820 seconds. SPP method  need run 24 iterations in 0.3274 seconds.

Figure 11 shows the results for computing Pareto H-eigenvalues of $\mathcal{A}$ from {\bf \it Example 6}, and the starting point is $x_0=[0.1846;0.8337;0.1696;0.9532;0.7225]$. In this case, all of the three methods can find the same Pareto H-eigenvalue 25.6537. SPG2 method need run 14 iterations in  0.4368 seconds. And SPG1 method need run 17 iterations in 0.2964 seconds while SPP method  need run 27 iterations in 0.3276 seconds.

We also used 100 random starting guesses for finding Pareto H-eigenvalue,  to test SPG1, SPG2, and SPP, respectively.  For each set of experiments, the same set of random starts was used. We listed the median number of iterations until convergence, and the average run time in the 100 experiments in Table 5. The computed Pareto H-eigenvalues were listed in Fig.8, Fig.10, Fig.12 for Example 4,5,6, respectively. As we can see from the Table 5, SPP is slightly slower than SPG method. SPP method need much more iterations in general. The number of iterations of SPG2 is the least.  But, SPG1 is faster than SPG2 for the test problems.

\begin{center}
\small{Table 5. Comparison results for 100 random test on computing Pareto H-eigenvalues of $\mathcal{A}$ from {\bf \it Examples 4-6}.}\\
\begin{tabular}{|c|c|c|c|c|c|c|c|c|}
\hline \multicolumn{ 1}{|c|}{Algorithm}  &  \multicolumn{ 2}{|c|}{SPG1} &  \multicolumn{ 2}{|c|}{SPG2} &  \multicolumn{ 2}{|c|}{SPP} \\
\cline{2-7}
\multicolumn{ 1}{|c|}{Example} & Its. &  Time & Its. &  Time & Its. &  Time\\
\hline
Ex. 4 &  22.94 & 0.3861 &  22.51 & 0.5934&  39.21 & 0.5203 \\
Ex. 5 &  21.67 & 0.3844 &  13.08 & 0.4345&  24.94 & 0.3354 \\
Ex. 6 &  17.99 & 0.3151 &  11.09 & 0.3175&  23.98 & 0.3260 \\
\hline
\end{tabular}
\end{center}



\section{Conclusion}

 In this paper, two monotone ascent spectral projected gradient algorithms were investigated for the tensor eigenvalue complementarity problem (TEiCP).  We also presented a shifted scaling-and-projection algorithm, which is a great improvement of the original SPA method \cite{LHQ2015-COAP1}. Numerical experiments show that spectral projected gradient methods are efficient and competitive to the shifted projected power method.

\section*{Acknowledgements}

 This work was supported in part by the National Natural Science Foundation of China (No.61262026, 11571905, 11501100), NCET Programm of the Ministry of Education (NCET 13-0738), JGZX programm of Jiangxi Province (20112BCB23027), Natural Science Foundation of Jiangxi Province (20132BAB201026), science and technology programm of Jiangxi Education Committee (LDJH12088), Program for Innovative Research Team in University of Henan Province (14IRTSTHN023).

\end{document}